\documentclass[12pt]{article}

\usepackage{amsmath, amsthm, amscd, amsfonts, amssymb}

\def\unu{\mbox{$1\!\!\,\rule{0,1mm}{2,3mm}\,$}}
\def\calz{{\mathcal{Z}}}

\def\({\left(}
\def\){\right)}
\def\pf{\n{\bf Proof.} }
\def\vsp{\vspace*{1,5mm}\\ }
\def\vspp{\vspace*{2mm}\\ }

\def\bk{\bigskip }
\def\mk{\medskip }
\def\sk{\smallskip }

\def\n{\noindent }
\def\dd{\displaystyle}
\def\hf{{\hfill\rule{1,5mm}{1,5mm}}}
\def\D{{\Delta}}

\def\barr{\begin{array}}
\def\earr{\end{array}}

\def\bit{\begin{itemize}}

\def\eit{\end{itemize}}
\def\D{{\Delta}}

\newtheorem{theorem}{Theorem}[section]

\newtheorem{corollary}[theorem]{Corollary}

\newtheorem{lemma}[theorem]{Lemma}

\newtheorem{remark}[theorem]{Remark}

\def\1{^{-1}}
\def\vsp{\vspace*{2mm}\\ }

\def\calf{{\mathcal{F}}}

\def\rr{{\mathbb{R}}}

\def\9{{\infty}}
\def\G{{\Gamma}}
\def\lbb{{\lambda}}

\def\b{{\beta}}

\def\g{{\gamma}}
\def\wt{\widetilde}
\def\ov{\overline}
\def\vf{{\varphi}}
\def\oo{{\omega}}
\def\ooo{{\Omega}}
\def\pp{{\partial}}
\def\D{{\Delta}}
\def\vp{{\varepsilon}}
\def\barr{\begin{array}}
\def\earr{\end{array}}

\def\dd{\displaystyle}
\def\bk{\bigskip }
\def\sk{\smallskip}
\def\n{\noindent }

\def\pas{\mathbb{P}\mbox{-a.s.}}

\def\vsp{\vspace*{2mm}\\ }
\def\ff{\forall }
\def\({\left(}
\def\){\right)}
\def\<{\left<}
\def\>{\right>}

\def\ns{Navier--Stokes }

\title{Global  solutions to random 3D vorticity equations for small initial data}
\author{Viorel Barbu\footnote{Octav Mayer Institute of
the Romanian Academy and Al.I. Cuza University, Ia\c si, Romania. Email: vb41@uaic.ro} \and Michael R\"ockner\footnote{Fakult\"at f\"ur Mathematik, Universit\"at Bielefeld, D-33501 Bielefeld, Germany. Email: roeckner@math.uni-bielefeld.de}}
\date{}

\begin{document}
\maketitle
\begin{abstract}
\n One proves the existence and uniqueness in $(L^p(\rr^3))^3$,
$\frac32<p<2$, of a    global mild solution to  random vorticity
equations associated to stochastic $3D$ Navier-Stokes equations
with linear  multiplicative Gaussian noise of convolution type,
for sufficiently small initial vorticity. This resembles some
earlier deterministic results of T. Kato \cite{6} and are obtained
by treating  the equation in vorticity form and reducing the
latter   to a random nonlinear parabolic equation. The solution
has  maximal regularity in the spatial variables and is weakly
continuous in $(L^3\cap L^{\frac{3p}{4p-6}})^3$ with respect to
the time variable.
Furthermore, we obtain the pathwise continuous dependence of solutions with respect to the initial data.\sk
\end{abstract}

\section{Introduction}\label{s1}
Consider the stochastic $3D$ \ns\ equation  \newpage

\begin{equation}\label{e1.1}
 \barr{ll}
dX{-}\D X\,dt+(X\cdot\nabla)X\,dt=\!\sum\limits^N_{i=1}
(B_i(X)+\lbb_iX)d\b_i(t){+}\nabla \pi\,dt\\
\hfill\mbox{on }(0,\9)\times\rr^3,\\ \nabla\cdot X=0\
\hfill\mbox{on }(0,\9)\times\rr^3,\vsp X(0)=x\ \hfill\mbox{in }
(L^p(\rr^3))^3, \earr\!\!\!\!\!\end{equation}where $\lbb_i\in\rr$,
$x:\ooo\to\rr^3$ is a random variable. Here $\pi$ denotes the
pressure and $\{\b_i\}^N_{i=1}$ is a system of independent
Brownian motions on a probability space $(\ooo,\calf,\mathbb{P})$
with normal filtration $(\calf_t)_{t\ge0}$, $x$ is
$\calf_0$--measurable, and $B_i$ are the convolution operators
\begin{equation}\label{e1.2}
B_i(X)(\xi)=\int_{\rr^3}h_i(\xi-\bar\xi)X(\bar\xi)d\bar\xi
=(h_i*X)(\xi),\ \xi\in\rr^3,
\end{equation}where $h_i\in L^1(\rr^3),\ i=1,2,...,N,$  and $\D$ is the Laplacian on $(L^2(\rr^3))^3$.

It is not known whether \eqref{e1.1} has a probabilistically
strong solution in the mild sense for all time. Therefore, we
shall rewrite \eqref{e1.1} in   vorticity form and then transform
it into a random PDE, which we shall prove, has a global in time
solution for $\mathbb{P}$-a.e. fixed $\oo\in\ooo$, provided the
initial condition is small enough.

Consider the vorticity field
\begin{equation}\label{e1.3}
U=\nabla\times X={\rm curl}\ X
 \end{equation}
and apply the curl operator to equation \eqref{e1.1}. We obtain (see e.g. \cite{2}, \cite{3}) the transport vorticity equation
\begin{equation}\label{e1.4}
 \barr{l}
dU-\D U\,dt+((X\cdot\nabla)U-(U\cdot\nabla)X)dt
=\sum\limits^N_{i=1}(h_i*U+\lbb_iU)d\b_i\vspp\hfill\mbox{ in
}(0,\9)\times\rr^3,\vspp U(0,\xi)=U_0(\xi)=({\rm curl}\,x)(\xi),\
\xi\in\rr^3.\earr \end{equation} The vorticity $U$ is related to
the velocity $X$ by the equation
\begin{equation}\label{e1.5}
X(t,\xi)=K(U(t))(\xi),\ t\in(0,\9),\ \xi\in\rr^3,\end{equation}
where $K$ is the Biot--Savart integral operator
\begin{equation}\label{e1.6}
K(u)(\xi)= -\frac1{4\pi}\int_{\rr^3}
\frac{\xi-\bar\xi}{|\xi-\bar\xi|^3}\times u(\bar\xi)d\tilde\xi,\
\xi\in\rr^3.\end{equation} Then one can rewrite the vorticity
equation \eqref{e1.4} as
\begin{equation}\label{e1.7}
\barr{l} dU-\D U\,dt+((K(U)\cdot\nabla)U-(U\cdot\nabla)K(U))dt
\vsp\qquad\quad=\sum\limits^N_{i=1}(h_i*U+\lbb_iU)d\b_i\
\hfill\mbox{ in }(0,\9)\times\rr^3,\vsp
  U(0,\xi)=U_0(\xi),\
\xi\in\rr^3.\earr\end{equation} Equivalently,
\begin{equation}\label{e1.7a}
\barr{ll} U(t)=\!\!\! &\dd e^{t\D}U_0-\!\! \int^t_0 \!e^{(t-s)\D}
((K(U(s)){\cdot}\nabla) U(s){-} (U(s){\cdot}\nabla)K(U(s)))ds\vsp
&+\dd\int^t_0\sum^N_{i=1}e^{-(t-s)\D}(h_i*U(s))+
\lbb_i(U(s))d\b_i(s),\ t\ge0.\earr\end{equation}

Now, we consider the transformation
\begin{equation}\label{e1.8}
U(t)=\Gamma(t)y(t),\ t\in[0,\9),\end{equation} where
$\Gamma(t):(L^2(\rr^3))^3\to(L^2(\rr^3))^3$ is the linear
continuous operator defined by the equations
\begin{equation}
d\Gamma(t) = \sum^N_{i=1} (B_i+\lbb_iI)\Gamma(t)d\b_i(t),\ t\ge0,\
\ \ \ \Gamma(0)=I,\label{e1.9}\end{equation} where (see
\eqref{e1.2}) \begin{equation} B_i u =h_i*u,\ \ff
u\in(L^p(\rr^3))^3,\ i=1,...,N,\ p\in(1,\9).\label{e1.10}
\end{equation}
We also set
\begin{equation}\label{e1.10a}
\wt B_i= B_i+\lbb_i I,\ i=1,...,N,\end{equation} where $I$ is the
identity operator.

Since $B_iB_j=B_jB_i$, equation \eqref{e1.9} has a solution $\Gamma$ and can be equivalently
expressed as (see~\cite{4}, Section 7.4)
\begin{equation}\label{e1.11}
\Gamma(t)=\prod^N_{i=1}\exp\left(\b_i(t)\wt B_i-\frac t2\ \wt
B^2_i\right),\ t\ge0.\end{equation}  Here \eqref{e1.9} is meant in
the sense that, for every $z_0\in (L^2(\rr^3))^3$, the continuous
$(\calf_t)$--adapted $(L^2(\rr^3))^3$--valued process
$z(t):=\G(t)z_0$, $t\ge0$, solves the following SDE on
$H:=(L^2(\rr^3))^3$,
$$dz(t)=\sum^N_{i=1}\wt B_iz(t)d\beta_i(t),\ \ z(0)=z_0,$$
where $H$ is equipped with the usual scalar product
$\left<,\cdot,\right>$.

Applying the It\^o formula in \eqref{e1.7} (the justification for
this is as in \cite{1a}), we obtain for $y$ the random
differential equation
\begin{equation}\label{e1.12}
\hspace*{-2mm}\barr{l}
\dd\frac{dy}{dt}\,(t)-\Gamma\1(t)\D(\Gamma(t)y(t))+
\Gamma\1(t)(K(\Gamma(t)y(t))\cdot\nabla)(\Gamma(t)y(t))
\vsp\qquad\quad \dd-(\Gamma(t)y(t)\cdot\nabla)(K(\Gamma(t)y(t)))
=0,\ \ t\in[0,\9), \vspp y(0)=U_0.\earr\end{equation} Taking
into account that, for all $i$, $B_i\D=\D B_i$ on $H^2(\rr^3)$, it
follows by \eqref{e1.9}, \eqref{e1.11} that
$\D\Gamma(t)=\Gamma(t)\D$ on $H^2(\rr^3)$, $\ff t\ge0$.

 In
what follows, equation \eqref{e1.12} will be taken in the
following mild sense
\begin{equation}\label{e1.14}
  \dd y(t)=e^{t\D}U_0
\dd+\int^t_0e^{(t-s)\D}\G\1(s)
 M(\G(s)y(s))ds,\ \
t\in[0,\9),
\end{equation}
where
\begin{equation}\label{e1.15}
(e^{t\D}u)(\xi)=\frac1{(4\pi t)^{\frac32}}\int_{\rr^3}\exp
\left(-\frac{|\xi-\bar\xi|^2}{4t} \right)u(\bar\xi)d\bar\xi,\
t\in[0,\9),\ \xi\in\rr^3,\end{equation} and $M$ is defined by
\begin{equation}\label{e1.16}
 \dd M(u)=- [(K( u)\cdot\nabla)(u)
-( u\cdot\nabla)(K(u))],\ t\in[0,\9). \end{equation}
We note that $U(t)= \Gamma(t) y(t)$ is the solution to the equation
\begin{equation}\label{e1.17a}
U(t) = e^{t\Delta} \Gamma(t) U_0 + \int^{t}_0 e^{(t-s)\Delta}
\Gamma(t)\Gamma^{-1}(s) M(U(s))ds ,\end{equation}
which may be viewed as the random version of the stochastic vorticity
equation \eqref{e1.7a}.

Our aim here and the principal contribution of this work is to
show that, for every $\vp\in(0,1)$, there exists
$\ooo_\vp\in\calf$   such that
$\mathbb{P}(\ooo_\vp)\ge1-\vp$ and, for all $\oo\in\ooo_\vp$, we
have the existence and uniqueness of a solution (in the mild
sense) for \eqref{e1.14} if the vorticity of $x$, i.e., $U_0={\rm
curl}\,x$, is $\pas$ sufficiently small in a sense to be made
precise   in Theorem \ref{t1.1} below.  We recall that, for a
deterministic Navier--Stokes equation, such a result was first
established by T.~Kato \cite{6} (see also T. Kato and H.~Fujita
\cite{6a})  and extended later to more general initial data by
Y.~Giga and T. Miyakawa \cite{5}, M. Taylor \cite{9}, H.~Koch and
D.~Tataru \cite{7}. However, the standard approach \cite{6},
\cite{6a} cannot be applied in the present case for one reason:
the nonlinear inertial term $(X\cdot\nabla)X$ cannot be
conveniently estimated in the space
$C_b([0,\9);L^p(\ooo\times\rr^d))$ and similarly for the
nonlinearity arising in \eqref{e1.7}. As regards the stochastic
$3D$ Navier-Stokes equations, to best of our knowledge all global
existence results were limited to martingale solutions. Since the fundamental work \cite{Flandoli}, the literature on (global) martingale solutions for stochastic $3D$-Navier-Stokes equations has grown enormously. We refer, e.g., to \cite{5aa}, \cite{Debusche}, \cite{Flandoli-1},  \cite{Goldys}, \cite{10a}, and the references therein.

In the following, we denote by $L^p$, $1\le p\le\9$, the space
$(L^p(\rr^3))^3$ with the norm $|\cdot|_p,$ by $W^{1,p}$ the
corresponding Sobolev space and by $C_b([0,\9);L^p)$ the space of
all bounded and continuous functions $u:[0,\9)\to L^p$ with the
sup norm. We  also set $D_i=\frac\pp{\pp\xi_i},$ $i=1,2,3$, and
denote by $\nabla\cdot u$
  the divergence of $u$, while
$$((u\cdot\nabla)v)_j=u_i D_iv_j,\ j=1,2,3,\ u=\{u_i\}^3_{i=1},\
v=\{v_j\}^3_{j=1}.$$As usual $$q'=\frac q{q-1}\mbox{\ \ for
}q\in(1,\9).$$
We set  for $p\in\(\frac32,3\)$
\begin{equation}\label{e1.18}
 \eta(t)=  \|\G(t) \|_{L(L^p,L^p)}
 \|\G(t)\|_{L(L^{\frac{3p}{3-p}},L^{\frac{3p}{3-p}})}
 \|\G\1(t)\|_{L^q,L^q)} ,\ t\ge0,
\end{equation}
 where  for $q\in (1,\9)$,
$\|\cdot\|_{L(L^q,L^q)}$ is the norm of the space $L(L^q,L^q)$ of
linear continuous operators on $L^q$.

For $p\in[1,\9)$, we denote by $\mathcal{Z}_p$ the space of all
functions $y:(0,\9)\times\rr^3\to\rr^3$ such that
\begin{equation}\label{e1.18a}
\barr{rcl}
t^{1-\frac3{2p}}y&\in&C_b([0,\9);L^{p}),\\
t^{\frac32\(1-\frac1p\)}D_iy&\in&C_b([0,\9);L^{p}),\ i=1,2,3.\earr
\end{equation}
The space $\mathcal{Z}_p$ is endowed   with the norm
\begin{equation}\label{e1.19}
 \hspace*{-2mm}\|y\|_{p,\9}= \dd\sup
\Big\{t^{1-\frac3{2p}}|y(t)|_{p}+
t^{\frac32\(1-\frac1p\)}|D_iy(t)|_{p};\, t\in(0,\9),\, i=1,2,3\Big\}.
\end{equation}
In the following, we take $\lbb_i\in\rr$  such that
\begin{equation}\label{e1.19a}
|\lbb_i|>(\sqrt{12}+3)|h_i|_1,\ \ \ff i=1,2,...,N.
\end{equation}
 We note that
$$\|B_i\|_{L(L^q,L^q)}\le|h_i|_{L^1},\ \ff i=1,...,N.$$
Theorem \ref{t1.1} is the main result.

 \begin{theorem}\label{t1.1} Let $p,q\in(1,\9)$   such that
\begin{equation}\label{e1.20}
\frac32<p<2,\ \ \frac1q=\frac2p-\frac13\,\cdot
 \end{equation}
Let $\ooo_0=\left\{\sup_{t\ge0}\eta(t)<\9\right\}$ and consider
\eqref{e1.14} for fixed $\oo\in\ooo_0$. Set $\G(t):=\G(t)(\oo),$
$\eta(t):=\eta(t,\oo)$. Then $\mathbb{P}(\ooo_0)=1$ and there is a
positive constant $C^*$ independent of $\oo\in\ooo_0$ such
that,~if $U_0\in L^{\frac32}$ is such that
\begin{equation}\label{e1.21}
\sup_{t\ge0} \eta(t) |U_0|_{\frac32} \le C^*,
\end{equation} then the random  equation \eqref{e1.14}
  has a unique solution $y\in\mathcal{Z}_p$
which satisfies
\begin{equation}
M(\G(t)y) \in L^1(0,T;L^q).  \label{e1.24}
\end{equation}
Moreover,  for each     $\vf\in L^{3}\cap L^{q'}$, the function $t\to
\int_{\rr^3}y(t,\xi)\vf(\xi)d\xi$ is continuous   on $[0,\9)$.
 The map $U_0\to y$ is   Lipschitz form $L^{\frac32}$ to~$\mathcal{Z}_p$.

 In particular, the random vorticity equation \eqref{e1.17a} has a unique solution $U$ such that $\G\1U\in \calz_p$.
\end{theorem}

\begin{remark}\label{r1.2}\rm Concerning condition \eqref{e1.21},
we note that an elementary calculation shows that
$$\eta(t)\le\prod^N_{i=1}\exp(3|\b_i(t)|
(|h_i|_1|+|\lbb_i|)-t\alpha_i),\
t\in[0,\9),$$ where $\alpha_i:=\frac12\
\lbb^2_i-\frac32\,(|h_i|^2_1+2|\lbb_i|\,|h_i|_1),$ which is
strictly positive by \eqref{e1.19a}.
\end{remark}

By the law of the iterated logarithm, it follows that
$$\sup_{t\ge0}\eta(t)<\9,\ \ \mathbb{P}\mbox{-a.e.},$$hence for $\ooo_r:=\{\sup_{t\ge0}\eta(t)\le r\}$ we have $\mathbb{P}(\ooo^c_r)\to0$ as $r\to\9.$

But, taking into account that, for each $r>0$ and all $\nu>0$, $i=1,...,N,$ we have (see Lemma 3.4 in \cite{1aa})
$$\mathbb{P}\left[\sup_{t\ge0}\left\{\exp(\b_i(t)-\nu t )\right\}\ge r\right]=r^{-2\nu},$$
 and, more explicitly,  we get that
$$
\mathbb{P}(\ooo^c_r)\le2N r^{-\frac{N\alpha}{\g^2}},\ \ \ff r>0,$$
where $\alpha=\min_{1\le i\le N}\alpha_i$,
$\g=3\max\{(|h_i|_1+|\lbb_i|);\ i\le N\}.$ Therefore, if $\oo\in\ooo_r$ and $U_0=U_0(\oo)$ satisfies
  \begin{equation}\label{e1.31}
 |U_0|_{\frac32}\le\frac{C^*}{r},\end{equation}
 then condition
\eqref{e1.21} holds.  It is trivial to
define such an $\calf$-measurable function $U_0:\ooo\to L^{\frac32}$,
for instance, $$U_0:=\sum^\9_{n=1}\frac{C^*}n\
\unu_{\{n-1\le\sup_{t\ge0}\eta(t)<n\}}u_0,$$for some $u_0\in
L^{\frac32}.$ But, of course, $U_0$  is not $\calf_0$-measurable
and so the process $U(t)$, $t\ge0$, given by Theorem \ref{t1.1},
is not $(\calf_t)_{t\ge0}$-adapted. Therefore, $U=\G(t)y$ is not a
solution to the stochastic vorticity equation \eqref{e1.7a}.
However, it can be viewed as a generalized solution to
\eqref{e1.7a}.

It should also be mentioned  that assumption \eqref{e1.19a} is not
necessary for existence of a solution to equation \eqref{e1.14},
but only to make sure that    condition \eqref{e1.21} is not void.

\section{Proof of Theorem \ref{t1.1}}
\setcounter{equation}{0}

To begin with, we note below in Lemma \ref{l2.1} a few immediate
properties of the operator $\Gamma$ defined in
\eqref{e1.9}--\eqref{e1.11}.
\begin {lemma}\label{l2.1} We have
\begin{equation}
|\Gamma(t)z|_q+|\Gamma\1(t)z|_q \le  C_t|z|_q,\ t\in[0,\9),\ \ff
z\in L^q,\ \ff q\in[1,\9),\label{e2.1} \end{equation} and
\begin{equation}\label{e2.2}
|\nabla(\G(t)z)|_q\le\|\G(t)\|_{L(L^q,L^q)}|\nabla z|_q,\mbox{ for all }z\in W^{1,q}(\rr^3). \end{equation}
\end{lemma}

\n{\bf Proof.} By \eqref{e1.2}, \eqref{e1.10} and by the Young
inequality, we see that
\begin{equation}
|B_i(u)|_q \le (|h_i|_1+|\lbb_i|)|u|_q,\ \ff u\in L^q,\
i=1,...,N.\label{e2.4}
\end{equation}

\n Recalling \eqref{e1.11}, we see by \eqref{e2.4} that
\eqref{e2.1}, \eqref{e2.2} hold,  as claimed.\hf

\begin{lemma}\label{l2.2} Let $\frac1q=\frac1{r_1}+\frac1{r_2}$,
$\frac32<r_1<\9,$ $r^*_1=\frac{3r_1}{3+r_1},$ $1<q<\9.$
 Then, for some $C>0$ independent of $\oo$,
\begin{equation}\label{e2.6}
\hspace*{-6mm}\barr{r}
 |M(\G(t)z)|_q\le C \|\G(t)\|_{L(L^{r_1},L^{r_1})}
\|\G(t)\|_{L(L^{r_2},L^{r_2})}
(|z|_{r_1}|z|_{r_2}+|z|_{r^*_1}|\nabla z|_{r_2}),\vsp
t\in[0,\9),\earr\!\!\!\!\!\end{equation} for all $z\in L^{r_1}\cap
L^{r_2}\cap L^{r^*_1}$ with $\nabla z\in L^{r_2}$.
\end{lemma}

\n{\bf Proof.} We have by \eqref{e1.16} and \eqref{e2.1}
\begin{equation}\label{e2.7}
|M(\G(t)z)|_q\le |(K(\G(t)z)\cdot\nabla)(\G(t)z)|_q
+|(\G(t)z\cdot\nabla)K(\G(t)z)|_q.
\end{equation}On the other hand, by \eqref{e2.2} and the H\"older inequality we have
\begin{equation}\label{e2.8}
\hspace*{-6mm}\barr{ll} |(K(\G(t)z)\cdot\nabla)(\G(t)z)|_q
\!\!\!&\le |K(\G(t)z)|_{r_1} |\nabla(\G(t)z)|_{r_2}\vsp &\le
\|\G(t)\|_{L(L^{r_1},L^{r_1})} \|\G(t)\|_{L(L^{r_2},L^{r_2})}
|K(z)|_{r_1}|\nabla z|_{r_2}.\earr\!\!\!\!\!\!
\end{equation}Now, we recall the classical estimate for Riesz
potentials (see \cite{8}, p.~119)
$$\left|\int_{\rr^3}\frac{f(\bar\xi)}
{|\xi-\bar\xi|^2}\,d\xi\right|_\beta \le C|f|_\alpha,\ \ff f\in
L^\alpha,$$where $\frac1\beta=\frac1\alpha-\frac13,$
$\alpha\in(1,3).$ By virtue of \eqref{e1.6}, this yields
\begin{equation}\label{e2.9}
|K(u)|_\beta\le C|u|_\alpha;\ \ff u\in L^\alpha,\
\frac1\beta=\frac1\alpha-\frac13,\end{equation}and so, for
$\beta=r_1$, $\alpha=\frac{3r_1}{3+r_1}=r^*_1$, we get by \eqref{e2.2} and
\eqref{e2.8} the estimate
\begin{equation}\label{e2.10}
|(K(\G(t)z)\cdot\nabla)(\G(t)z)|_q\le
C\|\G(t)\|_{L(L^{r_1},L^{r_1})}
\|\G(t)\|_{L(L^{r_2},L^{r_2})}|z|_{r^*_1}|\nabla z|_{r_2}.
\end{equation}
(Here and everywhere in
the following, $|\nabla z|_p$ means $\sup\{|D_iz|_p\,;\
i=1,2,3\}$.)
Taking into account that, by the Calderon--Zygmund inequality (see \cite{5a}, Theo\-rem 1),
\begin{equation}\label{e2.10a}
|\nabla K(z)|_{\tilde p}\le C|z|_{\tilde p},\ \ff z\in L^{\tilde
p},\ 1\le\tilde p<\9,\end{equation} one obtains that
  \begin{equation}\label{e2.11}
    |(\G(t)z\cdot\nabla)(K(t)\G(t)z)|_q
  \le C\|\G(t)\|_{L(L^{r_1},L^{r_1})}
\|\G(t)\|_{L(L^{r_2},L^{r_2})}|z|_{r_1}|z|_{r_2}.
\end{equation}
Substituting \eqref{e2.10}, \eqref{e2.11} in \eqref{e2.7}, one
obtains estimate \eqref{e2.6}, as claimed.\hf

\begin{lemma}\label{l2.3} Let $r_1=\frac{3r_2}{3-r_2},\
\frac32<r_2<3,$ $q=\frac{3r_1}{r_1+6}.$ Then, we have, for some
$C>0$ independent of $\oo$,
\begin{equation}\label{e2.12}
\hspace*{-3mm}|M(\G(t)z)|_q\le C\|\G(t)\|_{L(L^{r_1},L^{r_1})}
\|\G(t)\|_{L(L^{r_2},L^{r_2})}|z|_{r_2}|\nabla z|_{r_2},\, \ff
z\in W^{1,r_2}.\end{equation}
\end{lemma}

\pf We have by the Sobolev--Gagliardo-Nirenberg inequality (see,
e.g., \cite{2a}, p.~278)
$$|z|_{r_1}\le C|\nabla z|_{r_2},\ \ff z\in W^{1,r_2}(\rr^3).$$
Substituting in \eqref{e2.6} and taking into account that
$r^*_1=r_2$, we obtain \eqref{e2.12}, as claimed.\hf\mk

In the following, we fix $p=r_2$, $r_1$ and $q$ as in Lemma
\ref{l2.3}, \eqref{e2.12}, that is,
\begin{equation}\label{e2.13}
\frac32<p<2,\ r_1=\frac{3p}{3-p},\ \frac1q=\frac2p-\frac13.
\end{equation}
We write equation \eqref{e1.14} as
\begin{equation}\label{e2.14}
y(t)=G(y)(t)=e^{t\D} U_0+F(y)(t),\
t\in[0,\9),\end{equation}where
\begin{equation}\label{e2.15}
F(z)(t) =  \int^t_0  e^{(t-s)\D}\G\1(s)
M(\G(s)z(s))ds,\ t\in[0,\9).\end{equation} By \eqref{e1.15}, we have  for
$1<\tilde q\le\tilde p<\9$  the estimates
\begin{eqnarray}
|e^{t\D}u|_{\tilde p}&\le&C t^{-\frac32\left(\frac1{\tilde
q}-\frac1{\tilde p}\right)}|u|_{\tilde q},\ \ u\in L^{\wt q},\label{e2.16}\\[2mm]
|D_je^{t\D}u|_{\tilde p}&\le&C t^{-\frac32 \left(\frac1{\tilde
q}-\frac1{\tilde p}\right)-\frac12}|u|_{\tilde q},\
j=1,2,3.\label{e2.17}
\end{eqnarray}
(Everywhere in the following, we shall denote by $C$ several
positive constants independent of $\oo$ and $t\ge0$.)

  We apply   \eqref{e2.16}, with $\tilde q=q,\ \tilde p=p$. By
\eqref{e2.12}--\eqref{e2.15}, we obtain that
\begin{equation}\label{e2.18}
 \barr{ll}
 |F(z(t))|_{p}
 \dd \!\!\!&\le  C \!\dd\int^t_0 \!(t{-}s)^{-\frac12
 \(\frac3p-1\)}|\G\1(s)M(\G(s)z(s))|_q
  ds\vsp
 &\le  \dd
 C\!\!\int^t_0\!\!(t{-}s)^{-\frac12\(\frac3p-1\)}
 \|\G(s)
 \|_{L(L^{\frac{3p}{3-p}},L^{\frac{3p}{3-p}})}\vsp
 &\qquad\quad\|\G(s)\|_{L(L^p,L^p)}
 \|\G\1(s)\|_{L(L^q,L^q)}
 |z(s)|_p|\nabla z(s)|_p\,ds. \earr
\end{equation}
Similarly,  we obtain by \eqref{e2.17} that
\begin{equation}
\barr{l}
|D_jF(z(t))|_{p}\label{e2.19}
 \dd\dd\le C\! \int^t_0\!  (t{-}s)^{-\frac3{2p}}
\|\G(s)\|_{L(L^{\frac{3p}{3-p}},L^{\frac{3p}{3-p}})}\vsp
\qquad\quad\|\G(s)\|_{L(L^p,L^p)}
\|\G\1(s)\|_{L(L^q,L^q)} |z(s)|_p|\nabla z(s)|_p\, ds,\ j=1,2,3.
 \earr\end{equation}

We consider the Banach space $\mathcal{Z}_p$ defined by
\eqref{e1.18a}, that is,
\begin{equation}\label{e2.22}
\barr{r} \calz_p=\big\{y;\ t^{1-\frac3{2p}}y\in C_b([0,\9);L^p),\
t^{\frac32\(1-\frac1p\)}D_jy\in C_b([0,\9);L^p),\\
j=1,2,3\big\},\earr\end{equation} with the   norm
\begin{equation}\label{e2.23}
\|z\|_{p,\9}=\|z\|
=\!\dd\sup_{t>0}\left\{\(t^{1-\frac3{2p}}|z(t)|_p
+t^{\frac32\(1-\frac1p\)}|D_iz(t)|_p\),\,
  i=1,2,3\right\}.\end{equation} We note
that
\begin{equation}\label{e2.24}
|z(t)|_{p} |\nabla z(t)|_p\le C t^{-\frac52+\frac3p}\|z\|^2,\ \ff
z\in\calz_p,\ t\in(0,\9).\end{equation}
By   \eqref{e2.18}  and
\eqref{e2.23} we see that, for $z\in\calz_p$, we have
\begin{equation}\label{e2.25}
\hspace*{-3mm}\barr{ll}
 |F(z(t))|_{p}  \!\!\!
 &\le
\dd\int^t_0(t-s)^{-\frac12\(\frac3p-1\)}
\|\G(s)\|_{L(L^{\frac{3p}{3-p}},L^{\frac{3p}{3-p}})}
\|\G(s)\|_{L(L^p,L^p)}
 \vsp
 &\qquad\qquad\|\G\1(s)\|_{L(L^q,L^q)}
 |z(s)|_p|\nabla z(s)|_p\,ds\vsp
 &\le
C\dd\int^t_0(t-s)^{-\frac12\(\frac3p-1\)} |s|^{-\frac52+\frac3p}
\dd\|\G(s)\|_{L(L^{\frac{3p}{3-p}},L^{\frac{3p}{3-p}})}\vsp
&\qquad\qquad\|\G(s)\|_{L(L^p,L^p)}
 \|\G\1(s)\|_{L(L^q,L^q)}\,ds\|z\|^2 \vsp
   &\le C
t^{\frac3{2p}-1}\dd\sup_{0\le s\le t}\eta(s)
\dd\int^1_0(1-s)^{-\frac12\(\frac3p-1\)} s^{-\frac52+\frac3p}
ds\,\|z\|^2,\ \ff t>0 ,\earr\hspace*{-15mm}
\end{equation} where $\eta$ is given by \eqref{e1.18}.
This yields
\begin{equation}\label{e2.26}
t^{1-\frac3{2p}} |F(z(t))|_{p}\le \!C \sup_{0\le s\le t}
\{\eta(s)\}B\!\(\frac32\!\!\(\!\frac2p-1\)\!,
\frac32\(\!1-\frac1p\)\!\!\)\!\|z\|^2,\, \ff t{>}0,\ \!\!\!
\end{equation}where $B$ is the classical beta function (which is
finite by virtue of \eqref{e1.20}).

Similarly, by  \eqref{e2.17} and \eqref{e2.24}, we have, for
$j=1,2,3,$
\begin{equation}\label{e2.27}
\hspace*{-4mm}\barr{ll}|D_jF(z)(t)|_{p}
\!\!\!&\le C\dd\int^t_0(t-s)^{-\frac3{2p}}
 s^{-\frac52+\frac3p}
\|\G(s)\|_{L(L^{\frac{3p}{3-p}},L^{\frac{3p}{3-p}})}\vsp
&\qquad\quad
\|\G(s)\|_{L(L^p,L^p)} \|\G\1(s)\|_{L(L^q,L^q)}ds\|z\|^2\vsp
 & \le
 C \dd\sup_{0\le s\le t}\{\eta(s)\}t^{-\frac32\(1-\frac1p\)}
 B\(3\(\dd\frac1p-\frac12\),1-\frac3{2p}\)
 \|z\|^2,\ t>0. \earr\hspace*{-15mm}
 \end{equation} Hence,
 \begin{equation}\label{e2.28}
 t^{\frac32\(1-\frac1p\)}|D_jF(z(t))|_p\le C
 \sup_{0\le s\le t}\eta(s)\|z\|^2,\ \ff z\in\calz_p,\ t>0,\ j=1,2,3.\end{equation}
 By \eqref{e2.16}--\eqref{e2.17}, we have
 $$\barr{rcll}
 |e^{t\D} U_0|_p&\le&Ct^{\frac3{2p}-1}| U_0|_{\frac32},&t>0,\vspp
 |D_je^{t\D} U_0|_p&\le&Ct^{\frac3{2p}-\frac32}| U_0|_{\frac32},&t>0,\ j=1,2,3.\earr$$
 Therefore, by \eqref{e2.23}  we get
\begin{equation}\label{e2.29}
\|e^{t\D} U_0\| \le C | U_0|_{\frac32}.
\end{equation}
By \eqref{e2.23}, \eqref{e2.26}, \eqref{e2.28}, \eqref{e2.29}, we
get,
\begin{equation}\label{e2.30}
\|G(z)\|\le C_1\( | U_0|_{\frac32} +\sup_{t\ge0} \eta(t) \|z\|^2\)
,\ \ff z\in\calz_p,\end{equation}where $C_1>0$ is independent of
$\oo$.

We set \begin{eqnarray}
\eta_\infty&=&\sup_{t\ge0}\eta(t),\label{e2.32}\end{eqnarray} and
so \eqref{e2.30} yields
\begin{equation}\label{e2.33}
\|G(z)\|\le C_1(|U_0|_{\frac32}+\eta_\9\|z\|^2),\ \ff
z\in\calz_p.\end{equation}
We set
$$\Sigma=\{z\in\calz_p;\ \|z\|\le R^*\}$$and note that, by
\eqref{e2.33}, it follows that $G(\Sigma)\subset\Sigma$ if
\begin{equation}\label{e2.34}
|U_0|_{\frac32}\eta_\9\le(4C_1^2)\1,\
\end{equation}
(so $U_0$ must depend on $\oo$) and $R^*=R^*(\oo)$ is given by
\begin{equation}\label{e2.35}R^*=2C_1|U_0|_{\frac32}.\end{equation}
(We recall that $C_1$ is independent  of $\oo$ and $U_0.$) Moreover, by \eqref{e1.16} and \eqref{e2.15},
we have, for all $z,\bar z\in\calz_p$,
$$\barr{ll}
G(z)(t)  -\ G(\bar z)(t) =-\dd\int^t_0e^{(t-s)\D}\G\1(s)
[(K\G(s)(z(s)-\bar z(s))\cdot\nabla)\G(s)z(s) \vsp  +\
(K(\G(s)\bar z(s))\cdot\nabla)\G(s)(z(s)-\bar
z(s))-\G(s)(z(s)-\bar z(s))\cdot\nabla K(\G(s)z(s))\vsp  -\
(\G(s)\bar z(s)\cdot\nabla)K(\G(s)(z(s)-\bar z(s)))]ds.\earr$$
Proceeding as above, we get, as in \eqref{e2.18}, \eqref{e2.25},
\eqref{e2.26}, that
\begin{equation}\label{e2.36a}
\barr{ll}
\hspace*{-10mm}|G(z)(t){-}G(\bar z)(t)|_p\vsp
 \le C
\!\dd\int^t_0\!(t{-}s)^{-\frac12\(\frac3p-1\)}
  \dd\|\G(s)\|_{L(L^{\frac{3p}{p-3}},L^{\frac{3p}{p-3}})} \|\G(s)\|_{L(L^p,L^p)}\vsp
 \qquad\quad
   \|\G\1(s)\| _{L(L^q,L^q)}(|z(s)
  -\bar z(s)|_p (|\nabla
z(s)|_p+|\nabla\bar z(s)|_p)   \vsp
 \qquad\quad+|\nabla z(s)-\nabla\bar
z(s)|_p(|z(s)|_p+|\bar z(s)|_p))ds\vsp
 \le C
t^{-\(1-\frac3{2p}\)}\dd\sup_{0\le s\le t}\eta(s)\|z-\bar
z\|(\|z\|+\|\bar z\|),\ \ff t>0,\earr\end{equation} and also
(see \eqref{e2.19}, \eqref{e2.27}, \eqref{e2.28})
$$|D_jG(z)(t)-D_jG(\bar z(t))|_p\le
C t^{-\frac32\(1-\frac1p\)}
\sup_{0\le s\le t} \eta(s) \|z-\bar z\|(\|z\|+\|\bar z\|),\
\ff t>0,$$for $j=1,2,3.$ Hence, by \eqref{e2.23} and
\eqref{e2.32}, we obtain that
\begin{equation}\label{e2.36}\|G(z)-G(\bar
z)\|\le C_2\eta_\9R^*\|z-\bar z\|,\ \ff z,\bar
z\in\Sigma,\end{equation} where $C_2$ is independent of $\oo$.

Then, by \eqref{e2.35}, \eqref{e2.36}, it follows that, if
\eqref{e2.34}   and
\begin{equation}\label{e2.37}2C_1C_2\eta_\9|U_0|_{\frac32}<1,\end{equation}
hold,  then  the   operator $G$ is a contraction on   $\Sigma$
and so there is a unique solution $U\in\Sigma$ to 
\eqref{e1.14} provided \eqref{e1.21} holds with $C^*<(2C_1C_2)\1.$

Now, as seen earlier, by \eqref{e2.12}, \eqref{e1.14} and
\eqref{e2.24} we have
\begin{equation}\label{e2.38}
\hspace*{-5mm}\barr{ll}
 \dd|M(\G(t)y(t))|_q \!\!\!&\dd\le
C \|\G(t)\|_{L(L^{\frac{3p}{3-p}},L^{\frac{3p}{3-p}})}
  \|\G(t)\|_{L(L^p,L^p)}
  |y(t)|_{p}|\nabla y(t)|_{p}  \vsp
&\dd\le C\|\G(t)\|_{L(L^{\frac{3p}{3-p}},L^{\frac{3p}{3-p}})}
  \|\G(t)\|_{L(L^p,L^p)}
   t^{-\frac52+\frac3p}\|y\|^2,\,
\ff t>0.\earr\hspace*{-4mm}
\end{equation}
 On the other hand,  we have for all $\vf\in
 C^\9_0(\rr^3),$
 \begin{equation}\label{e2.39}
 \barr{ll}
\dd\int_{\rr^3}y(t,\xi) \cdot \vf(\xi)d\xi
=\dd\int_{\rr^3}(e^{t\D}) U_0(\xi) \cdot \vf(\xi)d\xi \\ \qquad+
\dd\int^t_0 \int_{\rr^3}\G\1(s)
  M(\G(s)y(s)) \cdot e^{(t-s)\D}
\vf(\xi)d\xi\,ds.\earr\!\!\!\end{equation} Recalling that, for all
$1\le \tilde p<\9$, $|e^{t\D}\vf|_{\tilde p}\le|\vf|_{\tilde p}, $
it follows by   \eqref{e2.38}   that
\begin{equation}
\barr{ll} \dd\left|\dd\int^t_0\!\!\int_{\rr^3}\!\!\G\1(s)
  M(\G(s)y(s)){\cdot}e^{(t-s)\D}
\vf(\xi)d\xi\,ds\right|
\\
\qquad\qquad\dd\le  C\dd\sup_{0\le s\le t}\eta(s)\dd\int^t_0\!\!\!
s^{-\frac52+\frac3p}   ds \|y\|^2
\, |\vf|_{q'}\label{e2.41}\vsp \qquad\qquad\le\dd C\dd\sup_{0\le
s\le t}\eta(s)\ t^{\frac3p-\frac32} \|y\|^2\,|\vf|_{q'},\, \ff
t\in(0,\9).\earr\end{equation} We also have by \eqref{e2.29}
$$\left|\int_{\rr^3}e^{t\D}  U_0(\xi) \vf(\xi)d\xi\right|\le
C|U_0|_{\frac32}  |\vf|_{3},\ \ff t\in[0,\9).$$
 Combining the latter with \eqref{e2.39},
 \eqref{e2.41},
 we obtain that, for $T>0$,
$$
\left|\dd\int_{\rr^3}y(t,\xi){\cdot}\vf(\xi)d\xi\right| \le C
T^{\frac3p-\frac32} (|\vf|_{q'}+|\vf|_3),\ \ff\vf\in L^{q'}\cap
L^3,\ t\in[0,T].$$ Hence, by \eqref{e2.39} and since $t\to
e^{t\D}U_0$ is continuous on $L^{\frac32}$, the function $t\to y(t)$ is
$L^3\cap L^{q'}$ weakly continuous on $[0,\9)$, where $q'=\frac{3p}{4p-6}$.

If $U=y(t,U_0)\in\calz_p$ is the solution to \eqref{e1.14},
equivalently \eqref{e2.14}, we have for all $U_0,\ov U_0$
satisfying \eqref{e1.21} (see \eqref{e2.29} and \eqref{e2.36})
$$\barr{l} \|y(\cdot,U_0)-y(\cdot,\ov U_0)\|  \le
\|e^{t\D} (U_0-\ov U_0)\| +\|F(y(\cdot,U_0))-F(y(\cdot,\ov
U_0))\| \vspp
 \qquad\qquad\ \le C  |U_0-\ov U_0|_{\frac32} +
 \eta_\9R^*C_2\|y(\cdot,U_0)-y(\cdot,\ov U_0)\| .\earr$$
Recalling that by \eqref{e2.35} and \eqref{e2.37} we have
$R^*C_2\eta_\9<1,$ this yields
 $$\barr{ll}
 \|y(\cdot,U_0)-y(\cdot,\ov U_0)\|
 \le\dd\frac C{1-R^*C_2\eta_\9}\
 |U_0-\ov U_0|_{\frac32}
  \le C(\oo) |U_0-\ov U_0|_{\frac32} , \earr$$and so, the  map $y\to U(\cdot,U_0)$ is   Lipschitz from $L^{\frac32}$
 to $\calz_p$. This completes the proof of Theorem \ref{t1.1}. \hf\bk

It should be noted that, by \eqref{e2.34} and \eqref{e2.35}, we have by the Fernique
theorem
$$|U_0|_{\frac32},R^*\in\bigcap_{r\ge1} L^r(\ooo),$$
and so, taking into account that $y\in\Sigma$, we see by
\eqref{e2.22}, \eqref{e2.23} that
\begin{equation}\label{e2.39a}
\sup_{t\ge0}\left\{t^{1-\frac3{2p}}|y(t)|_p
+t^{\frac32\(1-\frac1p\)}|D_iy(t)|_p\right\}
\in\bigcap_{r\ge1}L^r(\ooo),\ i=1,2,3.\end{equation} We have,
therefore, the following completion of Theorem \ref{t1.1}.

\begin{corollary}\label{t2.3} Under the assumptions of
Theorem {\rm\ref{t1.1}}, the solution $y=y(t,\oo)$ to the equation
\eqref{e1.14} satisfies   \eqref{e2.39a}. The same result holds
for the solution \mbox{$U(t)=\Gamma(t) y(t)$} of the random
vorticity equation \eqref{e1.17a}.
\end{corollary}

\section{The random version of the $3D$ Navier-Stokes equation}
\setcounter{equation}{0}

We fix in \eqref{e1.1} the initial random variable $x$ by the
formula
\begin{equation}\label{e4.1}
x=K(U_0),\end{equation} where $U_0={\rm curl}\,x$, $U_0=U_0(\oo)$
satisfies condition \eqref{e1.21} for all $\oo\in\ooo_0$ (see
Remark \ref{r1.2}). If $y$ is the corresponding solution to
equation \eqref{e1.14} given by Theorem \ref{t1.1}, we define the
process $X$ by formula \eqref{e1.5}, that~is,
\begin{equation}\label{e4.2}
X(t)=K(U(t))=K(\G(t)y(t)),\ \ \ff t\in[0,\9).\end{equation} By
\eqref{e2.9}, where $U$ is the solution to the vorticity equation
\eqref{e1.1}, we have
\begin{equation}\label{e4.3}|X(t)|_{\frac{3p}{3-p}}\le C |U(t)|_{p},\ \ \ff t\in[0,\9). \end{equation} (Everywhere in the following, $C$ are positive
constants independent of $\oo\in\ooo$.)

On the other hand, by the Carlderon--Zygmund ine\-qua\-lity
\eqref{e2.10a}, we have
\begin{eqnarray}
|D_iX(t)|_{p}&\le& C| U(t)|_{p},\
i=1,2,3.\label{e4.4}\end{eqnarray} By  \eqref{e4.3} and by Theorem \ref{t1.1}, it follows that
\begin{equation}
t^{1-\frac3{2p}} X \in C_b([0,\9);L^{\frac{3p}{3-p}}), \label{e4.5}
\end{equation}
while, by \eqref{e4.4}, we have for $i=1,2,3$
\begin{equation}\label{4.6}
t^{\frac32\(1-\frac1p\)}D_iX\in C_b([0,\9);L^{p}).
 \end{equation}

Now, if in \eqref{e2.10a} we take $z=D_jX$, we get  that, for all
$i,j=1,2,3$,
$$t^{\frac32\(1-\frac1p\)}|D_iD_jX|_{p}\le Ct^{\frac32\(1-\frac1p\)}|D_iU|_{p}\le C,\ \ff
t\in[0,\9).$$This yields
\begin{eqnarray}
t^{\frac32\(1-\frac1p\)}D^2_{ij}X&\in&C_b([0,\9);L^{p}),\
i,j=1,2,3.  \label{e4.8a}\end{eqnarray} Moreover,   by
Corollary \ref{t2.3}, we also have
\begin{eqnarray}
t^{1-\frac3{2p}}X&\in&C_b([0,\9);L^r(\ooo;L^{\frac{3p}{3-p}})),\
\ff
r\ge1,\label{e4.8}\\[1mm]
t^{\frac32\(1-\frac1p\)}D_iX&\in&C_b([0,\9);L^r(\ooo;L^p)),\
i=1,2,3,\label{e4.9}\\[1mm]
t^{\frac32\(1-\frac1p\)}D_{ij}X&\in&C_b([0,\9);L^r(\ooo;L^p)),\
i,j=1,2,3.\label{e4.10}
\end{eqnarray}

Now, if in equation \eqref{e1.17a} one applies the Biot--Savart
operator $K$, we obtain for $X$ the equation
\begin{equation}
 \hspace*{-5mm}\barr{r}
X(t)=K(e^{t\D}\G(t){\rm curl}\,x)
+\!\!\dd\int^t_0\!\!K(e^{(t-s)\D} \G(t)\G\1(s) M({\rm
curl}\,X(s)))ds,\\ t\ge0,\earr\!\!\!\label{e3.11}
\end{equation}
where $M$ is given by \eqref{e1.16}. It should be noted that, by virtue of \eqref{e4.8a}-\eqref{e4.10}, the right hand side of \eqref{e3.11} is well defined.

Equation \eqref{e3.11} can be viewed as the random version of the
Navier-Stokes equation \eqref{e1.1}. However, since, as seen earlier,  $U_0$ is not
$\calf_0$-measurable, the processes $t\to y(t)$, $t\to U(t)$ are not
$(\calf_t)_{t\ge0}$-adapted, and so $X$ is not
$(\calf_t)_{t\ge0}$-adapted, too. Therefore, \eqref{e3.11} cannot be
transformed back into \eqref{e1.1}. By Theorem \ref{t1.1}, it
follows that

\begin{theorem}\label{t3.1}
Under assumptions \eqref{e1.21}, the random Navier-Stokes equation
\eqref{e3.11} has a unique solution $X$ satisfying
\eqref{e4.8a}-\eqref{e4.10}.
\end{theorem}

\begin{remark}\rm As easily seen from the proofs,  Theorem \ref{t1.1}  extends mutatis-mutandis to the noises
$\sum\limits^N_{i=1}\sigma_i(t,X)\dot\beta_i(t)$, where
$$\sigma_i(t,x)(\xi)=\int_{\rr^3}
h_i(t,\xi-\bar\xi)x(\bar\xi)d\bar\xi,\ \xi\in\rr^3,\ i=1,...,N,$$
where $t\to h_i(t,\xi)$ is continuous and $$|h_i(t)|_1\le C,\ \ \ff t\ge0,\ i=1,...,N.$$\end{remark}

\begin{remark}\rm The linear multiplicative case $B_i(X):=\alpha_iX$, $i=1,...,N,$ that is $h_i:=\delta$, where $\delta$ is the Dirac measure, can be approximated by taking
$h_i(\xi)=\frac1{\vp^d}\ \rho\(\frac\xi\vp\)$, where $\rho\in
C^\9_0(\rr^d)$, support $\rho\subset\{\xi;\ |\xi|_d\le1\}$,
$\int\rho(\xi)d\xi=1.$  \end{remark}

\bk\n{\bf Acknowledgements.} Financial support through SFB701 at Bielefeld University is gratefully acknowledged. V. Barbu was also partially supported by a CNCS UEFISCDI (Romania) grant, project PN-II-ID-PCE-2012-4-0156.

\end{document}